# Mathematics of Balanced Parentheses: The case of Ordered Motzkin Words

Gennady Eremin

ergenns@gmail.com

December 20, 2019

**Abstract.** It is known that the core of mathematics is natural numbers. And everything related to the natural number is interesting to mathematicians. In this paper, we draw parallels between natural numbers and elements of a non-numeric lexicographic sequence, Motzkin words (well-formed strings of parentheses and zeros). We will also talk about the decomposition of Motzkin words into pairs of parentheses (analogue of Prime Numbers). Finally, we will try to interest the reader in the elements of mathematical analysis on bracket expressions. The nesting procedure of parentheses is described by the author as the differentiation of the weight function. We will work a little with derivatives for bracket pairs and also give some differential equations.

## 1   Introduction

In discrete mathematics, bracket sets are often studied. Of particular interest is the ordering of parentheses words, the analysis of disjoint and nested words, and the establishment of distances between words (for example, see [BP14] and [GZ14]). The goal of this paper is to construct a lexicographic series from Motzkin words, which will be as close as possible to a series of natural numbers. This allows you to enter arithmetic operations and even the basics of mathematical analysis on items of a sequence (a kind of derivatives and differential equations). This work is a continuation of [Ere19].

First, let's analyze the sequence of natural numbers on a purely formal level.

1.1. **Natural numbers and bitwise arithmetic.** Often zero is considered a natural number. This suits us, because among Motzkin words there is an identical element − the word "0". Let's write the set of natural numbers this way:

$$\mathbb{N} = \{0, 1, 2, …, 9, 10, 11, …, 99, 100, 101, …, 999, 1000, 1001 …\}.$$

We note several properties of set; it is these properties that interest us first of all.

(i)   The set is ordered in ascending order of numbers. Elements are indexed from zero; the index of each element is equal to its value: $n_i = i$, $i \geq 0$. Thus, natural numbers are self-indexing.

(ii)  Since indexes are not repeated, the elements of $\mathbb{N}$ are unique. There are no repeats among natural numbers.

(iii) At first glance, natural numbers are listed in order of increasing code length. Integers are distributed over *ranges*. Single-digit numbers, 1-*range*, are listed first, followed by double-digit integers, 2-*range*, and so on. Let's call it a *primary order*.

(iv) In ranges, integers are sorted according to the following order of weight of the alphabet symbols: $0 < 1 < 2 < ... < 9$. The minimum weight is 0, the maximum is 9. It is logical to call such sorting an *internal order*.

(v) A natural number can begin with an arbitrary digit except 0. The only first number 0 is an exception to the rule. In other cases, the symbol with the minimum weight is not written at the beginning of natural numbers.

As for the latter property, at times for convenience when performing some operations, we temporarily write a zero (or several zeros) before the natural number. But such temporary *leading zeros* do not change the number itself.

We will implement the considered properties of natural numbers on the elements of the lexicographic sequence. It remains for us to describe simple bitwise arithmetic, the analogue of which is used in brackets.

Bitwise arithmetic is the preprocessing of operands before addition or subtraction. In particular, before summing multi-bit numbers, we can adjust the digits where there are no carries. For example, $27 + 59 = 77 + 9$ (preprocessed digits are shown in red). Of course, for natural numbers, such bitwise procedures are not essential, but such processing of bracket sets is significant.

1.2. **Unique Motzkin Words.** Any Motzkin word is made up of *bricks* (*atoms*) of three types: the left parenthesis, the right parenthesis, and zero. Balanced parentheses mean that each opening bracket has a corresponding closing bracket; the left parenthesis precedes the right one. Such *matched pairs* are properly nested. The matched pairs of parentheses are similar to *molecules* that are located among zeros to compose the Motzkin word. In this sense, paired parentheses resemble prime numbers into which natural numbers are decomposed.

In combinatorics, we usually count the number of elements of some set. The set of Motzkin words of length $n$ are enumerated by the Motzkin numbers $M_n$ (see [OEIS A001006](https://oeis.org/A001006)). For example, there are two Motzkin words of size 2, *2-word*: 00 and ().The three-character Motzkin word, *3-word*, can be obtained in four variants ($M_3 = 4$): 000, 0(), (0), ()0. The first two 3-words are *inherited* from the 2-words by adding leading zero, while the last 3-words (0) and ()0 are *unique*. Next, we have four inherited 4-words and five unique 4-words: (00), (0)0, (()), ()00, ()().

Obviously, among the *n*-words, $M_{n-1}$ are inherited, and the remaining $M_n - M_{n-1} = U_n$, $n > 1$, are unique. The numbers $U_n$ for $n = 1, 2, \ldots$ form the sequence

1, 1, 2, 5, 12, 30, 76, 196, 512, 1353, 3610, 9713, 26324, 71799 …

We need to order Motzkin words like natural numbers. But we are hindered by items with leading zeros, inherited words. These words do not satisfy the last property of natural numbers. Therefore, it was decided not to consider words with leading zeros (the exception is the initial word "0").

Finally, we need a total alphabetical order; this order is logical: 0 < ( < ). In the Motzkin word, any matched pair starts with the left parenthesis, so it's logical to take the weight of the left parenthesis less than the right one. The symbol "0" resembles zero in integers (it has a minimum weight and it is free).

As a result, we get the set of the *Unique Motzkin Words*, UMWs (see the first 800 items in [Ere19]):

$$\mathfrak{M} = \{0, (), (0), ()0, (00), (0)0, (()), ()00, ()(), (000), (00)0, (0()), (0)00 \ldots\}$$

In $\mathfrak{M}$, UMWs are indexed from zero; let's denote them like this: $w_0 \equiv 0$, $w_1 \equiv ()$, $w_2 \equiv (0)$, and so on. We identify $w_i$ ($i$ is a specific index) and the $i$-th item in $\mathfrak{M}$. In both sets $\mathfrak{M}$ and $\mathbb{N}$, initial items are identical, $w_0 = n_0$. These two elements have similar properties.

In $\mathfrak{M}$, all elements are distributed over ranges along the code length: $\mathfrak{M}_1 = \{0\}$, $\mathfrak{M}_2 = \{()\}$, $\mathfrak{M}_3 = \{(0), ()0\}$, and so on. The cardinality of the $n$-range is $\#\mathfrak{M}_n = U_n$, $n > 1$. In the $n$-range, $n > 1$, we have the minimum (shown in red) $\min \mathfrak{M}_n = (0^{n-2})$ and the maximum $\max \mathfrak{M}_n = ()^{\lfloor n/2 \rfloor}[0]$. The superscript indicates the repetition of zero or a matched pair. The maximum $n$-word ends with 0 if $n$ is odd.

Let's call the index of $x \in \mathfrak{M}$ *weight* and denote wt $x$. So, wt $(\min \mathfrak{M}_n) = M_{n-1}$, wt $(\max \mathfrak{M}_n) = M_n - 1$, wt $w_i = i$, and so on. In general, the weight of an item does not always match its index. For example, 5-word $w_{11} \equiv (0())$ is made up of two components $w_9 \equiv (000)$ and $w_3 \equiv ()0$. When we put $w_3$ inside $w_9$, the *nest-weight* of $w_3$ changes, namely wt' $w_3 = 2$. We will talk about nesting below.

## 2   Maths of prime words

A matched pair of parentheses and everything inside is called a *block*. A block that is not contained within another block is called a *prime word*. A prime word can have ending zeros. Let's solve the task of decomposition of any UMW into prime words. Below we describe *partial addition* $\oplus$ and *partial subtraction* $\ominus$.

Let the 9-word $w_{736} \equiv ()0(0())0$ be given. In $w_{736}$, there are two prime words: 9-word $()0^7 \equiv w_{708}$ and 6-word $(0())0 \equiv w_{28}$. At once we get the *weight expression* $708 + 28 = 736$. For Motzkin words we get the equality $w_{736} = w_{708} \oplus w_{28}$. Also we can write down corresponding subtraction operations: $w_{708} = w_{736} \ominus w_{28}$ and $w_{28} = w_{736} \ominus w_{708}$.

When performing arithmetic operations, zero in Motzkin words is processed like numerical zero. Let's write down it in the form of the following rules:

$$0 \oplus 0 = 0, \quad 0 \oplus ( = ( \oplus 0 = (, \quad 0 \oplus ) = ) \oplus 0 = ),$$
$$0 \ominus 0 = 0, \quad ( \ominus 0 = (, \quad ) \ominus 0 = ), \quad ( \ominus ( = 0, \quad ) \ominus ) = 0.$$

Obviously, for $x \in \mathfrak{M}$   $x \oplus 0 = 0 \oplus x = x$, $x \ominus 0 = x$, $x \ominus x = 0$. Operation $\oplus$ occur if the operands do not *intersect*. Let's check the weight of the 9-word $()0^7$. Since $()0^7 \oplus \max \mathfrak{M}_7 = \max \mathfrak{M}_9$, the following weight expression is true:

$$\text{wt } ()0^7 = \text{wt}(\max \mathfrak{M}_9 \ominus \max \mathfrak{M}_7) = (M_9 - 1) - (M_7 - 1) = 834 - 126 = 708.$$

The considered operations allow concatenating prime words. Above we met the 6-word (0())0; to receive such codes you need to be able to nest brackets.

## 3 Prime pairs

Let's call a prime word with a single pair of parentheses a *prime pair*. Using prime pairs we can compose any Motzkin word. Let's write the corresponding set:

$\mathscr{P}$ = { (), (0), ()0, (00), (0)0, ()00, (000), (00)0, (0)00, ()000, (0000), (000)0 …}

In $\mathscr{P} \subset \mathfrak{M}$, elements are indexed from one: $p_1$ = (), $p_2$ = (0), $p_{11}$ = (0000), and so on. Items are distributed over ranges; the cardinality of the *n*-range is $\#\mathscr{P}_n = n-1$. Let's denote a prime pair of size *n* with the right parenthesis in the *k*-th position as $p_{n,k} = (0^{n-k-1})0^{k-1}$, $n > k > 0$. The first *n*-word (shown in red) is $p_{n,1} = (0^{n-2})$, the last *n*-word is $p_{n,n-1} = ()0^{n-2}$. It is easy to find the index of any pair in $\mathscr{P}$: $p_{n,k} = p_i$, $i = k + 1 + 2 + … + (n-1) = k + (n-1)(n-2)/2$.

What is the index (or rather the weight) of $p_{n,k}$ in $\mathfrak{M}$? Looking through $\mathfrak{M}$, we get: wt $p_{n,1} = M_{n-1}$, wt $p_{n,n-1} = M_n - M_{n-2}$. Let's get the general formula starting with $p_{k+1,k}$, and then we will increase the word length by moving the left bracket (see Corollary 3.1 in [Ere19]):

$p_{k+1,k} = ()0^{k-1}$, wt $p_{k+1,k} = M_{k+1} - M_{k-1}$;
$p_{k+2,k} = (0)0^{k-1}$, wt $p_{k+2,k} = M_{k+1} - M_{k-1} + (M_{k+1} - M_k)$;
$p_{k+3,k} = (00)0^{k-1}$, wt $p_{k+3,k} = M_{k+1} - M_{k-1} + (M_{k+2} - M_k)$;
……..
$p_{n,k} = (0^{n-k-1})0^{k-1}$, wt $p_{n,k} = M_{k+1} - M_{k-1} + (M_{n-1} - M_k)$.

As a result, we get

(3.1) $\qquad$ wt $p_{n,k} = M_{n-1} + U_{k+1} - M_{k-1}$, $n > k > 0$.

In General,

(3.1a) $\qquad$ wt $p_{n,k} = M_{n-1} + \delta_k$, $\delta_k = U_{k+1} - M_{k-1}$, $n > k > 0$.

The numbers $\delta_k$ for $k = 1, 2 …$ form the sequence

$\qquad$ 0, 1, 3, 8, 21, 55, 145, 385, 1030, 2775, 7525, 20526, 56288 …

## 4 Nesting procedure

Prime pairs can be *nested*, and we must be able to calculate the weight of the resulting Motzkin words. The procedure for including a prime pair into a Motzkin word and the corresponding processing of indices resemble the author the *differentiation* of a function. Below are a few steps.

Let be a prime pair $p_{n,k}$, $n > k > 1$ (in the case $k = 1$, the pair cannot be nested anywhere), then

$$\text{wt } p_{n,k} = M_{n-1} + M_{k+1} - M_k - M_{k-1}.$$

Next, we increase both *n* and *k* by 1

$$\text{wt } p_{n+1,\,k+1} = M_n + M_{k+2} - M_{k+1} - M_k,$$

and calculate the change in weight:

$$\Delta \text{wt } p_{n,\,k} = \text{wt } p_{n+1,\,k+1} - \text{wt } p_{n,\,k} = U_n + M_{k+2} - 2M_{k+1} + M_{k-1}.$$

And finally, we reduce the obtained value by $M_k$ and call it the *derivative* of function wt at point $p_{n,\,k}$, that is,

(4.1) $\quad \text{wt}' p_{n,\,k} = \Delta \text{wt } p_{n,\,k} - M_k = U_n + M_{k+2} - 2M_{k+1} - U_k, \ n > k > 1.$

What is interesting about (4.1)? We constructed an expression that gives an increment in weight, a *nest-weight* with which $p_{n,\,k}$ is included in the Motzkin word. Note we are talking about the first level of nesting, since the first derivative is obtained.

Let's test (4.1) for the 6-word $w_{28} \equiv (0())0$ that we met before:

$$\text{wt } (0())0 = \text{wt } (000)0 + \text{wt}' ()00$$
$$= \text{wt } p_{6,2} + \text{wt}' p_{4,3}$$
$$= (M_5 + 1) + (U_4 + M_5 - 2M_4 - U_3)$$
$$= (21 + 1) + (5 + 21 - 2 \times 9 - 2) = 28.$$

It is not difficult to calculate the derivative for some values of $k$. For example, $\text{wt}' p_{n,\,2} = U_n$, $\text{wt}' p_{n,\,3} = U_n + 1$, and so on. In General,

(4.1a) $\quad \text{wt}' p_{n,\,k} = U_n + \delta'_k, \ \delta'_k = M_{k+2} - 2M_{k+1} - U_k, \ n > k > 1.$

The numbers $\delta'_k$ for $k = 2, 3 \ldots$ form the sequence

$$0, 1, 4, 13, 39, 113, 322, 910, 2562, 7203, 20251, 56980, 160524 \ldots$$

## 5  Differential equations

Let's denote $\text{wt}^{(0)} x = \text{wt } x$, then (4.1) can be written as a *first-order differential equation* of the form

(5.1) $\quad \text{wt}' p_{n,\,k} + \text{wt}^{(0)} p_{n,\,k} + M_k = \text{wt}^{(0)} p_{n+1,\,k+1}, \ n > k > 1.$

Let's call the resulting differential equation *three-in-one*. There is a general form of the three-in-one equation. Here is the differential equation of the $(s+1)$-th order:

(5.2) $\quad \text{wt}^{(s+1)} p_{n,\,k} + \text{wt}^{(s)} p_{n,\,k} + \text{wt}^{(s-1)} p_{n,\,k} = \text{wt}^{(s)} p_{n+1,\,k+1}, \ n > k > s+1.$

Let's check 12-word $w_{9763} \equiv ((00)0(0()))$:

$$\text{wt } ((00)0(0())) = \text{wt } (0^{10}) + \text{wt}' (00)0^7 + \text{wt}' (000)0 + \text{wt}'' ()00$$
$$= \text{wt } p_{12,1} + \text{wt}' p_{11,8} + \text{wt}' p_{6,2} + \text{wt}'' p_{4,3}.$$

We calculate the summands using the equalities (3.1a), (4.1a) and (5.2):

$\text{wt } p_{12,1} = M_{11} = 5798,$
$\text{wt}' p_{11,8} = U_{11} + \delta'_8 = 3610 + 322 = 3932,$
$\text{wt}' p_{6,2} = U_6 = 30,$

$$\text{wt}'' p_{4,3} = \text{wt}' p_{5,4} - \text{wt}' p_{4,3} - \text{wt}\, p_{4,3}$$
$$= (U_5 + \delta'_4) - (U_4 + \delta'_3) - (M_3 + \delta_3) = (12+4) - (5+1) - (4+3) = 3.$$

As a result we get the desired result

$$\text{wt}\, w_{9763} = 5798 + 3932 + 30 + 3 = 9763.$$

Using equation (5.2), it is easy to calculate the nest-weight of pairs. Below in Appendix we have given a table of the initial nest-weights of prime pairs from ten ranges. An inquisitive reader can independently obtain additional equations by analyzing the information of Appendix. For example, here is the obvious equality

$$\text{wt}^{(n-2)} p_{n,\,n-1} = n - 1.$$

Using the nest-weights we can easy decompose a Motzkin word into prime pairs (an analog of the factorization of natural numbers). In this case, the matched pair is characterized by three parameters: the positions of the left and right brackets, plus the depth of the nest. The outer brackets of prime words have zero nesting depth.

The reverse procedure can generate an arbitrary Motzkin word in accordance with the given parameters of prime pairs.

Gzhel State University, Moscow, 140155, Russia
http://www.en.art-gzhel.ru/

**Appendix**. Prime Pairs and their *nest*-weights (starting weights for ten ranges).

| №№ | n/k | $p_{n,k}$ | $M_k$ | wt | wt′ | wt″ | wt‴ | wt$^{iv}$ | wt$^{v}$ |
|---|---|---|---|---|---|---|---|---|---|
| 1 | 2/1 | ( ) | 1 | 1 | – | – | – | – | – |
| 2 | 3/1 | (0) | 1 | 2 | – | – | – | – | – |
| 3 | 3/2 | ( )0 | 2 | 3 | 2 | – | – | – | – |
| 4 | 4/1 | (00) | 1 | 4 | – | – | – | – | – |
| 5 | 4/2 | (0)0 | 2 | 5 | 5 | – | – | – | – |
| 6 | 4/3 | ( )00 | 4 | 7 | 6 | 3 | – | – | – |
| 7 | 5/1 | (000) | 1 | 9 | – | – | – | – | – |
| 8 | 5/2 | (00)0 | 2 | 10 | 12 | – | – | – | – |
| 9 | 5/3 | (0)00 | 4 | 12 | 13 | 9 | – | – | – |
| 10 | 5/4 | ( )000 | 9 | 17 | 16 | 10 | 4 | – | – |
| 11 | 6/1 | (0000) | 1 | 21 | – | – | – | – | – |
| 12 | 6/2 | (000)0 | 2 | 22 | 30 | – | – | – | – |
| 13 | 6/3 | (00)00 | 4 | 24 | 31 | 25 | – | – | – |
| 14 | 6/4 | (0)000 | 9 | 29 | 34 | 26 | 14 | – | – |
| 15 | 6/5 | ( )0000 | 21 | 42 | 43 | 30 | 15 | 5 | – |
| 16 | 7/1 | (00000) | 1 | 51 | – | – | – | – | – |
| 17 | 7/2 | (0000)0 | 2 | 52 | 76 | – | – | – | – |
| 18 | 7/3 | (000)00 | 4 | 54 | 77 | 69 | – | – | – |
| 19 | 7/4 | (00)000 | 9 | 59 | 80 | 70 | 44 | – | – |
| 20 | 7/5 | (0)0000 | 21 | 72 | 89 | 74 | 45 | 20 | – |
| 21 | 7/6 | ( )00000 | 51 | 106 | 115 | 88 | 50 | 21 | 6 |
| 22 | 8/1 | (000000) | 1 | 127 | – | – | – | – | – |
| 23 | 8/2 | (00000)0 | 2 | 128 | 196 | – | – | – | – |
| 24 | 8/3 | (0000)00 | 4 | 130 | 197 | 189 | – | – | – |
| 25 | 8/4 | (000)000 | 9 | 135 | 200 | 190 | 133 | – | – |
| 26 | 8/5 | (00)0000 | 21 | 148 | 209 | 194 | 134 | 70 | – |
| 27 | 8/6 | (0)00000 | 51 | 182 | 235 | 208 | 139 | 71 | 27 |
| 28 | 8/7 | ( )000000 | 127 | 272 | 309 | 253 | 159 | 77 | 28 |
| 29 | 9/1 | (0000000) | 1 | 323 | – | – | – | – | – |
| 30 | 9/2 | (000000)0 | 2 | 324 | 512 | – | – | – | – |
| 31 | 9/3 | (00000)00 | 4 | 326 | 513 | 518 | – | – | – |
| 32 | 9/4 | (0000)000 | 9 | 331 | 516 | 519 | 392 | – | – |
| 33 | 9/5 | (000)0000 | 21 | 344 | 525 | 523 | 393 | 230 | – |
| 34 | 9/6 | (00)00000 | 51 | 378 | 551 | 537 | 398 | 231 | 104 |
| 35 | 9/7 | (0)000000 | 127 | 468 | 625 | 582 | 418 | 237 | 105 |
| 36 | 9/8 | ( )0000000 | 323 | 708 | 834 | 721 | 489 | 264 | 112 |
| 37 | 10/1 | (00000000) | 1 | 835 | – | – | – | – | – |
| 38 | 10/2 | (0000000)0 | 2 | 836 | 1353 | – | – | – | – |
| 39 | 10/3 | (000000)00 | 4 | 838 | 1354 | 1422 | – | – | – |
| 40 | 10/4 | (00000)000 | 9 | 843 | 1357 | 1423 | 1140 | – | – |
| 41 | 10/5 | (0000)0000 | 21 | 856 | 1366 | 1427 | 1141 | 726 | – |
| 42 | 10/6 | (000)00000 | 51 | 890 | 1392 | 1441 | 1146 | 727 | 369 |
| 43 | 10/7 | (00)000000 | 127 | 980 | 1466 | 1486 | 1166 | 733 | 370 |
| 44 | 10/8 | (0)0000000 | 323 | 1220 | 1675 | 1625 | 1237 | 760 | 377 |
| 45 | 10/9 | ( )00000000 | 835 | 1865 | 2263 | 2044 | 1474 | 865 | 412 |
| 46 | 11/1 | (000000000) | 1 | 2188 | – | – | – | – | – |